\documentclass[12pt]{amsart}
\usepackage{latexsym}
\usepackage{amssymb}
\usepackage{amsmath}
\usepackage{amsfonts}
\usepackage{amscd}
\usepackage{graphicx}
\usepackage[all]{xy}
\usepackage{t1enc}

\newcommand{\n}{\mathbb{N}}

\newcommand{\h}{\mathbb{H}}
\newcommand{\f}{\mathbb{F}}
\theoremstyle{plain}
\newtheorem{lemme}{{Lemma}}[]
\newtheorem{proposition}{{Proposition}}[]

\newtheorem{definition}{{Definition}}[]
\newtheorem{corollaire}{{Corollary}}[]

\newtheorem{example}{{Example}}[]

\newtheorem{note}{{Notes}}[]
\newenvironment{prof}[1][Proof]{\textbf{#1.} }{\ \rule{0.5em}{0.5em}}
\theoremstyle{plain}

\begin{document}
\title{ CATEGORY OF FUZZY HYPER BCK-ALGEBRAS}
\author{J.DONGHO*}
\address {****Department of Mathematics, University of  Yaounde, BP
812,  Cameroon} \email{\makeatletter josephdongho@yahoo.fr}
 \begin{abstract}
 In this paper we first define the category of fuzzy hyper
 BCK-algebras. After that we show that the category of hyper
 BCK-algebras has equalizers, coequalizers, products.
 It is a consequence that this category is complete and hence has
 pullbacks.
\end{abstract}
\maketitle
\section{Introduction}
The study of hyperstructure was initiated in 1934 by F. Marty at
8th congress of Scandinavian Mathematiciens. Y.B. Jun et al.
applied the hyperstructures to BCK-algebras, and introduces the
notion of hyper BCK-algebra. Now we follow [1,2,3,4] and introduce
the category of fuzzy hyperBCK-algebra and obtain some result, as
mentioned in the abstarct.

\section{Preliminaries}
We now review some basic definitions that are very useful in the
paper.
\begin{definition}$[\ref{3}]$ Let $H$ be an non empty set.\\
A hyperoperation $*$ on $H$ is a mapping of $H\times H$ family of
non-empty subsets of $H$ $\mathcal{P}^*(H)$\end{definition}
\begin{definition}Let $*$ be an hyperoperation on $H$ and $O$ a constant element of $H$
An hyperorder on $H$ is subset $<$ of $\mathcal{P}^*(H)\times
\mathcal{P}^*(H)$ define by:\\ for all $x,y\in H, x<y$ iff $O\in
x*y$ and for every $A,B\subseteq H, A<B$ iff $\forall a\in A,
\exists b\in B$ such that $a<b.$
\end{definition}
\begin{definition}If $*$ is hyperoperation on $H$. \\
For all $A,B \subseteq H, A*B:=\underset{a\in A, b\in B}\bigcup a*b$\end{definition}
\begin{definition}$[1]$ By hyper BCK-algebra we mean a non empty set
$H$ endowed with a hyper-operation $*$ and a constant $O$
satisfying the following axioms.
\begin{itemize}\item[(HK1)]$(x*z)*(y*z)<(x*y)$
\item[(HK2)] $(x*y)*z=(x*z)*y$ \item[(HK3)] $x*H<\{x\}$
\end{itemize}
\end{definition}
\begin{definition} A fuzzy hyper BCK-algebra is a pair $(\mathbf{H}; \mu_H)$
where $\mathbf{H}=(H;*;O)$ is hyper BCK-algebra and
$\mu_H:H\longrightarrow [0,1]$ is a map satisfy the following
property:
$$\inf(\mu_H(x*y))\geq \min(\mu_H(x),\mu_H(y)) $$ for all $x,y\in
H.$\end{definition}
\begin{example}$[\ref{5}]$
 Let $n\in\mathbb{N}^*.$ Define the hyperoperation $*$ on $H=[n,+\infty)$
as follows:
$$x*y=\left\{\begin{array}{lllc}[n,x]& \texttt{iff}\quad x<y\\
                                (n,y]& \texttt{iff}\quad x>y\neq n\\
                                \{x\}& \texttt{iff}\quad y=n\end{array}\right.$$

for all $x,y\in H$. To show that $(H,*,n)$ is hyper BCK-algebra,
it suffice to show axiom  $HK3.$ For all $x\in H,
x*H=\underset{t\in H}\bigcup x*t.$ For all $x\in H$ then
$x*x\subseteq x*H.$ And then $n\in [n, x]*\{x\}$
\end{example}
\section{The category of fuzzyhyper BCK-algebras}
\begin{lemme}Let  $(\mathbf{H}; \mu_H)$ be a fuzzy hyper BCK-algebra. \\
For all $x\in H,\mu_H(O)\geq\mu_H(x)$\end{lemme}
\begin{prof} For all $x\in H, x<x;$ then $O\in x*x$.\\
$$\begin{array}{ccll}O\in x*x& \texttt{imply}&
\mu_H(O)\geq\inf(\mu_H(x*y))\geq \min(\mu_H(x),\mu_H(y))\\
                            &\texttt{i.e}& \mu_H(O)\geq
                            \min(\mu_H(x),\mu_H(y))= \mu_H(x)\\
                            &\texttt{i.e}& \mu_H(O)\geq \mu_H(x).
\end{array}$$
\end{prof}
\begin{definition}Let  $(\mathbf{H};\mu_H)$ be a fuzzy hyperBCK-algebra. $\mu_H$ is called
a fuzzy map.
\end{definition}
\begin{lemme} Let $(\mathbf{H};\mu_H)$ a fuzzy hyper BCK-algebra.The following properties
are trues:\begin{itemize}\item[i)]If  for all $x,y\in H, x<y$
 imply $\mu_H(x)\leq \mu_H(y)$ \\then for all $x\in H, \mu_H(x)=\mu_H(O)$
 \item[ii)] If $\mu_H(O)=0$ then $\mu_H(x)=0$
\end{itemize}
\end{lemme}
\begin{prof}
\begin{itemize}
\item[i)] For all $x\in H$, $x*H<\{x\}$ then $x*O<x$.\\
$O<x\Rightarrow \mu_H(O)\leq \mu_H(x)$. \\Then $\mu_H(x)\leq
\mu_H(O)$ and $\mu_H(O)\leq \mu_H(x)$ for all $x\in H$. \\i.e
$\mu_H(x)=\mu_H(O)$ for all $x\in H$. \item[ii)]
$\mu_H(O)=O\Rightarrow \mu_H(O)\leq\mu_H(x),$ for all $x\in H.$
\\Then $\mu_H(x)=\mu_H(O)$ for all $x\in H$.
\end{itemize}\end{prof}
\begin{definition} Let $(\mathbf{H};\mu_H)$ and $(\mathbf{F},\mu_F)$ two fuzzy
hyperBCK-algebras. An homomorphism from $(\mathbf{H},\mu_H)$ to
$(\mathbf{F},\mu_F)$ is an homomorphism
$f:\mathbf{H}\longrightarrow\mathbf{F}$ of hyper BCK-algebra such
that for all $x\in H$, $\mu_F(f(x))\geq\mu_H(x)$
\end{definition}
\begin{proposition}Let $(\mathbf{H},\mu_H)$ an hyperBCK-algebra. Let $\mathbf{G},\mathbf{F}\subset H$
two hyperBCK-sualgebras of $\mathbf{H}$. If there exist $\alpha\in
]0,1[$ such that \\$\mu_H(G^*)\subset [0,\alpha[$ and
$\mu_H(F)\subseteq ]\alpha, 1]$. Then any homorphism of hyper
BCK-algebra $f:G\longrightarrow F$ is homomorphism of fuzzy
hyperBCK-algebra.
 \end{proposition}
 \begin{prof} Suppose that there is $\alpha\in ]0,1]$
 such that\\ $\mu_H(G^*)\subset [0,\alpha[$ and $\mu_H(F)\subseteq
 ]\alpha,1[.$\\
 Let $f:G\longrightarrow F$ an homomorphism of hyper
 BCK-algebra.\\
 For all $x\in G^*,f(x)\in F.$ And $\mu_F(f(x))>\alpha>\mu_H(x).$
 \\Then
 $\mu_F(f(x))>\mu_H(x)$ for all $x\in G^*$ \\$f(O)=O$ then
 $\mu_F(f(O))=\mu_F(O)=\mu_H(O)$\\ i.e $\mu_F(f(O))=\mu_H(O).$
 therefore, for all $x\in x\in G, \mu_F(x)\geq\mu_H(x)$
 \end{prof}
\begin{example}$[\ref{1}]$ Define the hyper operation $"*"$ on $H=[1; +\infty]$ as follow.
$$x*y=\left\{\begin{array}{ccccll}& [1,x]& \quad &\texttt{if}&\quad x\leq y\\
                                 &(1,y]&\quad\quad \quad&\texttt{if}&\quad x>y\neq 1\\
                              &\{x\}&\quad &\texttt{if}&\quad y=1\end{array}\right.$$
For all $x,y\in H, (\mathbf{H},*,1)$ is hyperBCK-algebra. Define
the fuzzy structure $\mu_H$ on $H$ by:
$$\begin{array}{lllcc}\mu_H:&H&\longrightarrow& [0,1]&\\
                          &x&\mapsto&\frac{1}{x}&\end{array}$$
We show that $(\mathbf{H},\mu_H)$ is a fuzzy hyper BCK-algebra.\\
Let $x,y\in H.$ \begin{enumerate}\item[(i)] If $x\leq y$, then
$x*y=[1,x]$; i.e for all $t\in x*y, 1\leq t\leq x\leq y$ \\and so
$\frac{1}{y}\leq \frac{1}{x}\leq \frac{1}{t}.$ So, $\mu_H(t)\geq
\frac{1}{y}=\min\{\frac{1}{y},\frac{1}{x}\}=\min\{\mu_H(x),\mu_H(y)\}$.\\
Then $\inf\{x*y\}\geq \min\{\mu_H(x),\mu_H(y)\}$ \item[(ii)] If
$x>y\neq 1$ then $x*y=(1,y]$. For all $t\in H\cap x*y,
\frac{1}{x}\leq \frac{1}{y}\leq \frac{1}{t}\leq 1.$ therefore,
$\mu_H(t)=\frac{1}{t}\geq \frac{1}{x}=\min\{\mu_H(x),\mu_H(y)\}$
for all $t\in x*y$. Then
$\in\{\mu_H(x*y)\}\geq\min\{\mu_H(x),\mu_H(y)\}.$ \item[(iii)] If
$y=1, x*y=\{x\}$, hence
$\mu_H(x*y)=\{\mu_H(x)\}=\{\frac{1}{x}\}$.\\ $y=1$imply $y\leq x$
and $\frac{1}{x}\leq \frac{1}{y}$ for all $x\in H$; i.e; \\
$\min\{\mu_H(x),\mu_H(y)\}=\frac{1}{x}$. Then
$\mu_H(x*y)=\{\frac{1}{x}\}.$ \\Thus
$\inf\{\mu_H(x*y)\}=\frac{1}{x}\geq \min\{(\mu_H(x),\mu_H(y))\}$
\end{enumerate}
\end{example}
\begin{proposition}The fuzzy hyperBCK-algebras and homomorphisms of fuzzy hyperBCK-algebras form a category.
\end{proposition}

\begin{prof} The proof is easy. \end{prof}

\begin{note} In the following we let $\mathcal{H}$ the category of hyperBCK-algebras; $\mathbb{F}_{\mathcal{H}}$ the
category of fuzzy hyperBCK-algebras; $\mathbb{H}$ the fuzzy hyper
BCK-algebra $(\mathbf{H},\mu_H)$
\end{note}

For any fuzzy hyper BCK-algebra $ \h$, we associate for all
$\alpha\in [0,1]$ the set $H_{\alpha}:=\{ x\in H,\mu_H(x)\geq
\alpha\}$
\begin{lemme} Let $\h$ a fuzzy hyper BCK-algebra. For all $\alpha\in [0,1], O\in H_{\alpha}$ and for all
$x,y\in H, x*y\subseteq H_{\alpha}$
 \end{lemme}
\begin{prof} By lemma 1, for all $x\in H, \mu_H(x)\leq \mu_H(0).$\\
 Then for all $x\in H_{\alpha}, \mu_H(O)\geq \mu_H(x)>\alpha$ i.e
 $O\in H_{\alpha}$.\\
 Let $x,y\in H_{\alpha}$; \\for all $t\in x*y,$ $\mu_H(t)\geq \inf\{\mu_H(x*y)\geq\min\{\mu_H(x),\mu_H(y)\}\}\geq
 \alpha$\\
 then $t\in H_{\alpha}$. therefore, $x*y\subseteq H_{\alpha}$

 \end{prof}

\begin{definition}Let $(H,*,O)$ be an hyper BCK-algebra. An hyper BCK-subalgebra of $H$
is a non empty subset $S$ of $H$ such that  $O\in S$ and $S$ is
hyper BCK-algebra with respect to the hyper operation $"*"$ on $H$

\end{definition}
\begin{proposition}Let $(H,*,O)$ be an hyper BCK-algebra. A non empty subset $S$ of $H$ is hyper BCK-subalgebra of $H$ iff for all
$x,y\in S, x*y\in S$ \end{proposition}
\begin{prof} The proof is easy.\end{prof}
\begin{definition} A fuzzy hyper BCK-subalgebra of $\h$ is an hyper BCK-subalgebra $S$ of $\mathbf{H}$
with the restriction $\mu_S$ of $\mu_H$ on $S.$
\end{definition}
\begin{proposition} For all $\alpha\in [0,1],$ $(H_{\alpha}, \mu_H)$ is fuzzy hyper BCK-subalgebra of $\h$
\end{proposition}
\begin{prof}By lemma 3, $H_{\alpha}$ is hyper BCK-subalgebra of $\mathbf{H}$ \\ and
$\inf\{ \mu_H(x*y)\}\geq \min\{\mu_H(x),\mu_H(y)\}$\end{prof}
\begin{definition} Let $\h$ by an fuzzy hyper BCK-algebra.
The fuzzy-hyperBCK-subalgebra $\mathbf{H}_{\alpha}:=(H_{\alpha};
\mu_H)$ is calling hyper $\alpha$-cut of $\h$ \end{definition}
\begin{proposition} Let $\h$ be fuzzy hyper BCK-algebra. A hyper BCK-subalgebra $S$ of $\mathbf{H}$ is fuzzy
hyper BCK-subalgebra iff $S$ is hyper $\alpha$-cut of $H.$
\end{proposition}
\begin{prof} By prosition 4, any hyper $\alpha$-cut is fuzzy hyper BCK-subalgebra. \\ Conversely, let $S$
be fuzzy hyper BCK-subalgebra of $\h$. Then $\mu_H(S)$ is subset
of $[0,1].$\\ If $0\in \mu_H(S),$ then $S=H_0=\h.$\\
If $0<\inf(\mu_H(S)),$ then $S=H_{\inf(\mu_H(S))}.$
\end{prof}
\begin{proposition} Let $\h$ and $\f$ be two fuzzy hyper BCK algebras. An $\mathcal{H}$-morphism
$f:H\longrightarrow F$ is $\f_{\mathcal{H}}$-morphism iff for all
$\alpha\in [0,1], f(H_{\alpha})\subseteq F_{\alpha}.$
 \end{proposition}

\begin{prof} Suppose that $f(H_{\alpha})\subseteq F_{\alpha}$
for all $\alpha\in [0,1]$ Let $x\in [0, 1]$ we need $\mu_H(x)\leq
\mu_F(f(x))$. Let $\alpha=\mu_H(x);\quad x\in H_{\alpha}$ and
$f(x)\in f(H_{\alpha})\subseteq F_{\alpha}$. Then
$\mu_F(f(x))>\alpha=\mu_H(x).$ whence for all $x\in H,
\mu_F(f(x))\geq \mu_H(x)$.\\ Conversely, suppose that
$f:\h\longrightarrow \f$ is $\f_{\mathcal{H}}$-morphism.\\
For all $x\in H_{\alpha}$ for some $\alpha \in [0, 1]$,
$\mu_F(f(x))\geq\mu_H(x)\geq \alpha$i.e; $f(x)\in [0, 1]$. Then
$f(H_{\alpha})\subseteq F_{\alpha}$ for all $\alpha\in [0, 1]$.

\end{prof}
\begin{proposition} A $\f_{\mathcal{H}}$-morphism $f:\h\longrightarrow\f$
is $\f_{\mathcal{H}}$-iso iff it is both $\mathcal{H}$-iso and
$\mu_H=\mu_F f.$

\end{proposition}
\begin{prof} Suppose that $f$ is $\mathcal{H}$-iso and $\mu_H=\mu_F f$.
there is $g\in Hom_{\mathcal{H}}(\mathbf{F},\mathbf{H})$; $g\circ
f=Id_H$ and $g\circ g =Id_F.$ \\ Then, for all $x\in F,
\mu_H(g(x))=\mu_F(f(g(x)))\mu_F(x).$ \\ And then, $g\in
Hom_{\f_{\mathcal{H}}}(\f,\h).$\\
Conversely, Suppose that $f$ is $\f_{\mathcal{H}}$-iso.\\ There is
$g\in Hom_{\f_{\mathcal{H}}}(\f,\h); g\circ f=Id_F$ and $f\circ
g=Id_H.$\\ Since $f\in Hom_{\f_{\mathcal{H}}}(\f,\h), \mu_H\leq
\mu_Ff.$\\ Since $g\in Hom_{\f_{\mathcal{H}}}(\h,\h), \mu_F\leq
\mu_Hg.$  $x\in H$ imply $f(x)\in F$. Then $\mu_F(f(x))\leq
\mu_H(g(f(x)))=\mu_H(x)$ i.e; $\mu_Ff\leq \mu_H.$\\ therefore,
$\mu_Ff=\mu_H$
\end{prof}
\begin{proposition} Let $f\in Hom_{\f_{\mathcal{H}}}(\f,\h).$\\
$f$ is $\h_\mathcal{H}$-mono iff $f$ is $\mathcal{H}$-mono
\end{proposition}
\begin{prof} Suppose that $f$ is $\h_\mathcal{H}$-mono.\\
 For all $h,g\in Hom_\mathcal{H}(\mathbf{K},\mathbf{H})$, such
 that $fh=fg,$ we define $\mu_K=\min (\mu_H(h(x);\mu_H(g(x))))$
 for all $x\in K.$
 \begin{enumerate} \item[a)] We show that $(K;\mu_K)$ is fuzzy hyper BCK-algebra.
 $$\begin{array}{lllccccllll} \inf(\mu_K(x*y))&=&\inf\{\mu_H(h(x*y); \mu_H(g(x*y)))\}\\
                                      &=&\inf\{ \mu_H(h(x)*h(y); \mu_H(g(x)*g(y)))\}\\
                                      &=&\min\{\inf\{\mu_H(h(x)*h(y)\};\inf\{\mu_H(g(x)*g(y)))\}\}\\
                                      &\geq &\min\{\min\{\mu_H(h(x); \mu_H(h(y)\};\min\{\mu_H(g(x);
                                      g(y))\}\}\\
                                      &\geq & \min\{\min\{\mu_K(x);
                                      \mu_K(y)\}\}\\
                                      &\geq &\min\{\mu_K(x); \mu_K(y)\}

 \end{array}$$ Then, for all $x,y\in K, \inf(\mu_K(x*y))\geq \min
 \{\mu_K(x),\mu_K(y)\}$. therefore, $(K; \mu_K)$ is fuzzy hyper
 BCK-algebra. \item[b)]We show that $h$ and $g$ are
 $\f_\mathcal{H}$-homomorphism.\\ For all  $x\in K,$ $\mu_K(x)=\min\{\mu_H(h(x),
 \mu_H(g(x)))\}.$ \\ Then $\mu_K(x)\leq \mu_H(g(x))$ and
 $\mu_K(x)\leq \mu_H(h(x))$.\\ therefore, $h$ and $g$ are
 $\f_\mathcal{}H$-morphism.\\ Since $f$ is $\f_\mathcal{H}$-mono
 and $h, g\in Hom_{\f_\mathcal{H}}(\f,\h),$ $fh=fg$ imply $f=g$\\
 Conversely, if $f$ is $\f_\mathcal{H}$-mono, it is
 $\mathcal{H}$-mono.
 \end{enumerate}
\end{prof}
\begin{lemme} The pair $\mathbb{O}=(\{O\}, \mu_o)$ where
$$\begin{array}{lcccll}\mu_o:&\{o\}&\longrightarrow & [0, 1]\\
&o&\mapsto & 0\end{array}$$ is fuzzy hyper BCK-algebra
\end{lemme}
\begin{prof} Easy\end{prof}
\begin{lemme} $\mathbb{O}$ is final objet of $\f_\mathcal{H}$\end{lemme}
\begin{proposition} The category $\f_\mathcal{H}$ has products.\end{proposition}
\begin{prof} Let $(\h_i;\mu_{H_i})_{i\in I}$ a family of fuzzy hyper
BCK-algebras.\\
Denote $\mathbf{H}=\underset{i\in I}{\prod}H_i$ the
$\mathcal{H}$-product of $(H_i)_{i\in I}$ with the projection
morphisms $p_i:\mathbf{H}\longrightarrow H_i$. Consider the
following map $\mu_H:H\longrightarrow[0,1]$ define
by:$$\mu_H(x)=\underset{i\in I}{\bigwedge}\mu_{H_i}p_i(x)$$ for
all $x\in H$
\begin{itemize}\item[a)]We show that  the pair $(\mathbf{H}; \mu_H)$ is fuzzy hyper
BCK-algebra.\\
For all $x,y\in H, p_i(x*y)=p_i(x)*p_i(y)$ for all $i\n I.$\\ Then
$$\begin{array}{lcllc}\inf(\mu_{H_i}p_i(x*y))&=&\inf(\mu_{H_i}(p_i(x)*p_i(y))\\
                                         &\geq& \min\{\mu_{H_i}(p_i(x));
                                         \mu_{H_i}(p_i(y))\}\end{array}$$

for all $i\in I$.\\ Then,
$$\begin{array}{lllcc}\inf(\underset{i\in
I}{\bigwedge}\mu_{H_i}p_i(x*y))&\geq &\underset{i\in
I}{\bigwedge}\inf\{\mu_{H_i}(p_i(x)*p_i(y)\}\\
                               &\geq &\underset{i\in
I}{\bigwedge}\min\{\mu_{H_i}(p_i(x));
                                         \mu_{H_i}(p_i(y))\}\\
                               &\geq & \min\{\underset{i\in
I}{\bigwedge}\mu_{H_i}p_i(x); \mu_{H_i}p_i(y)\}\\
                               &\geq & \min\{\mu_H(x),\mu_H(y)\}
\end{array}.$$
\item[b)] For all $i\in I, x\in H; \mu_{H_i}p_i(x)\geq
(\underset{i\in I}{\bigwedge}\mu_{H_i}p_i)(x)$.\\ Then each $p_i$
is $\f_\mathcal{H}$-morphism.\item[c)] If
$q_j:\f\longrightarrow\h_j$ is family of
$\f_\mathcal{H}$-morphism, there is unique $\mathcal{H}$-morphism
$\varphi:\mathbf{F}\longrightarrow\mathbf{H}$ such that the
following diagram commute. \[\xymatrix{H\ar[r]^{p_j}&H_j\\
                                     F\ar[u]^{\varphi}\ar[ur]_{q_j}&}\]
                                     i.e $p_j\varphi=q_j$ for all
                                     $j\in I$\\
For all $x\in F$, $\mu_F(x)\leq
\mu_{H_i}q_i(x)$.$$\begin{array}{lllcc}
\texttt{Then}&\mu_F(x)&\leq
\mu_{H_i}p_i\varphi(x) \quad \texttt{for all }\quad x\in F, i\in I\\
\texttt{i.e}&\mu_F(x)&\leq  \underset{i\in
I}{\bigwedge}\mu_{H_i}p_i\varphi (x) \\
&&\leq  (\underset{i\in I}{\bigwedge}\mu_{H_i}p_i)\varphi(x)\\
&&\leq \mu_H(\varphi(x))\quad \texttt{for all} \quad x\in F.\\
&\texttt{Then}&\mu_F\leq \mu_H\varphi\leq
\end{array}$$ Then $\varphi$ is $\f_\mathcal{H}$-morphism.

\end{itemize}

\end{prof}
\begin{proposition} $\f_\mathcal{H}$ have equalizers.\end{proposition}
\begin{prof}Let $f,g\in Hom_{\f_\mathcal{H}}(\h, \f), K:=\{ x\in H,f(x)=g(x)\}$.\\
It is prove in [1] that $K$ is hyper BCK-subalgebra of $H$. It is
clear that $(K,\mu_H)$ is fuzzy hyper BCK-algebra. Let
$i:K\longrightarrow H$ the inclusion map. $i\in
Hom_{\f_\mathcal{H}}(\mathbb{K}, \f).$ For all $x\in K,
fi(x)=f(x)=g(x)=gi(x).$\\
Let $h\in Hom_{\f_\mathcal{H}}(\mathbb{L}, \f)$ such that $fh=gh,$
for all $x\in L, f(h(x))=g(h(x))$. Then $Imh\subseteq L$. Define
$\delta:L\longrightarrow K$ by $\delta(x)=h(x)$ for all $x\in L.$
$\delta \in Hom_{\f_\mathcal{H}}(\mathbb{\h}, \f)$ and
$i\delta=h.$ So, the following diagram commute.
$$\xymatrix{\mathbb{K}\ar[r]^i&\h\ar@<1ex>[r]^f\ar[r]_g&\f\\
             \mathbb{L}\ar[u]^\delta\ar[ur]_h}$$
             Since $i$ is monic, $\delta$ is unique
             $\f_\mathcal{H}$-morphism such that the above diagram
             commute. \\ therefore, $\f_\mathcal{H}$ have
             equalizers.
\end{prof}
\begin{proposition} $\f_\mathcal{H}$ is complet.\end{proposition}
\begin{prof} By proposition 9, each family of objets of $\f_\mathcal{H}$ has
product.\\
 By proposition 10, each pair of parallel arrows has an equalizer. Then $\f_\mathcal{H}$ is complet.
 \end{prof}
 \begin{corollaire} $\f_\mathcal{H}$ has pulbacks \end{corollaire}
 \begin{prof}By propositions 9 and 10, $\f_\mathcal{H}$ has equalizers and products.
 therefore, $\f_\mathcal{H}$ has pulbacks.\end{prof}
\begin{proposition} $\mathcal{\f_H}$ have  coequalizers  \end{proposition}
\begin{prof}  Let $f,g\in
Hom_\mathcal{\f_H}(\h,\mathbb{K})$. Let $$\sum_{fg}=\{ \theta,
\theta\quad \texttt{regular congruence relation on}\quad
\mathbf{K}\quad \texttt{such that} \quad f(a)\theta g(a) \forall
a\in H\}$$ $\sum_{fg}\neq \phi$ because $K\times K\in \sum_{fg}$\\
Let $\rho =\underset{\theta\in\sum_{fg}}{\bigcap}\theta$. Then,
$\rho$ is regular congruence relation. \\ Define on  $K/\rho$ the
following hyper operation  $$[x]_\rho *[y]_\rho=[x*y]_{{\rho}}.$$
$(K/\rho;
*;[0]_\rho)$ is an objet of $\mathcal{H}$ (see $[\ref{1}]$).\\
 Define on $K/\rho$ the following map
 $$\begin{array}{lllcc}\mu_{K/\rho}:&K/\rho&\longrightarrow & [0,1]\\
                                        {}&\mu_{K/\rho}([x]_\rho)&\longmapsto&\underset{a\in
 [x]_\rho}\bigvee\mu_{K}(a)\end{array}$$
 \begin{itemize}
 \item[a)]We show that $(K/\rho, \mu_{K/\rho})$ is objet of $\f_\mathcal{H}$.\\
 If $x,y\in K$ such that
 $[x]_\rho=[y]_\rho$.\\ Then $$\underset{a\in[x]_\rho}\bigvee\mu_{K}(a)=\underset{a\in [y]_\theta}\bigvee\mu_{K}(a)
 $$\\
 $\forall x\in K,\,\mu_{K}(x)\leq \underset{a\in[x]_\rho}\bigvee
 \mu_{K}(a)$.\\
 Then $$\mu_{K}\leq \mu_{K/\rho}([x]_\rho)=\mu_{K/\rho}(\pi(x))$$
 Then, the canonical projection $\pi$ is an
 $\mathcal{\f_H}-morphism$\\Since for all $x\in H, f(x)\rho g(x)$,
 then $[f(x)]_\rho=[g(x)]_\rho$.\\therefore,  $(\pi\circ f)(x)=(\pi\circ
 g)(x).$\\Then, $\pi\circ f=\pi\circ g.$\item[b)]  Universal property of coequalizer.\\
 Let $\varphi:\mathbb{K}\longrightarrow \mathbb{L}$ and
 $\f_{\mathcal{H}}$-morphism such that $\varphi\circ f=\varphi\circ
 g$.\\ Define the following mapping.
 $$\begin{array}{lllcc}\psi:&K/\rho&\longrightarrow &L\\
                          {}&[x]_{\rho}&\longmapsto &\varphi(x) \end{array}$$
 \item[c)]We prove that $\psi$ is well define. \\If $[x]_\rho=[y]_\rho$
 then, for all $a\in H, \varphi(f(a))=\varphi(g(a))$ imply
 $f(a)R_{\varphi}g(a)$ because $R_{\varphi}$ is regular congruence
 on $K$. Then $R_\varphi\in \sum_{f,g}.$ The minimality of $\rho$
 on $\sum_{f,g}$ imply $\rho\subseteq R_\varphi$. \\therefore,
 $[x]_\rho=[y]_{\rho}$ imply $x\rho y$.\\ Then $x R_\varphi y$. i.e
 $\varphi(x)=\varphi(y)$\\And then, $\psi([x]_\rho)=\psi([y]_{\rho})$
 therefore, $\psi$ is well define.\\
 For all $x\in K, \mu_{L}(\psi(\pi)(x))=\mu_{L}(\varphi(x))\geq \mu_{K}(x), \forall x\in K$
then for all $a\in [x]_\rho$. $\mu_{L}(\varphi
(a))\geq\mu_{K}(a).$ By the minimality of $\rho$,
$[a]_\rho=[x]_\rho$ imply $a\rho x$ then $a R_\varphi x$ i.e
 $\varphi(a)=\varphi(x)$ (because $\rho\subseteq R_\varphi$). Then
 $$\underset{a\in[x]_\rho}\bigvee\tilde{L}(\varphi(a))=\underset{a\in[x]_\rho}{\bigvee}\tilde{L}(\varphi(x))
 =\tilde{L}(\varphi(x))$$ therefore
 $$\tilde{L}(\varphi(x))\geq\underset{a\in[x]_\rho}{\bigvee}(a)=\tilde{K}/\rho([x]_\rho)$$
  i.e $\tilde{L}(\psi([x]_\rho))\geq \tilde{k}/\rho([x]_\rho) \forall x\in H$. It is clean that
  $$\psi(\pi(x))=\psi([x]_\rho)=\varphi(x),\forall x\in H$$
  i.e $\psi\circ \pi=\varphi $\\ This prove the commutativity of the following diagram:
$$\xymatrix{H\ar@<2pt>[r]^f\ar@<-2pt>[r]_g_{.}&K\ar[r]^{\pi}\ar[dr]_\varphi&K/\rho\ar[d]^\psi\\&&L}$$
   The unicity of $\psi$ is thus to the fat that $\pi$ is
  epimorphism.\\
  Then, $\mathcal{F}_\mathcal{H}$ have coequalizer.\end{itemize} \end{prof}
\section*{ACKNOWLEDGEMENTS}

\end{document}